\documentclass{amsart}

\usepackage{amssymb}

\newtheorem{theorem}{Theorem}[section]
\newtheorem{lemma}[theorem]{Lemma}

\theoremstyle{definition}

\newtheorem{corollary}[theorem]{Corollary}
\newtheorem{proposition}[theorem]{Proposition}

\theoremstyle{remark}

\numberwithin{equation}{section}

 \usepackage{amssymb,latexsym,amsmath,amscd,mathrsfs,yfonts,hyperref}

\newcommand{\ra}{\rightarrow}

\DeclareMathOperator{\SO}{SO}

\DeclareMathOperator{\diag}{diag}

\DeclareMathOperator{\SL}{SL}

\DeclareMathOperator{\rank}{rank}

\title[]{An example of a non-amenable dynamical system which is boundary amenable} 

\begin{document}

\author{Jacopo Bassi}
\address{
Department of Mathematics, University of Tor Vergata, Via della Ricerca Scientifica 1, 00133 Roma, Italy}
\email{bssjcp01@uniroma2.it}

\author{Florin R\u adulescu*}
\address{Department of Mathematics, University of Tor Vergata, Via della Ricerca Scientifica 1, 00133 Roma, Italy}
\email{radulesc@mat.uniroma2.it}
\thanks{\\
*Florin R\u adulescu is a member of the Institute of Mathematics of the Romanian Academy}


\begin{abstract}
\noindent It is shown that the action of $\SL(3,\mathbb{Z})$ on the Stone-{\v C}ech boundary of $ \SL(3,\mathbb{Z}) / \SL(2,\mathbb{Z}) $ is amenable. This confirms a prediction by Bekka and Kalantar.





\end{abstract}
 
 \maketitle
\bibliographystyle{amsplain}

 
\section{A dynamical presentation of $\SL(3,\mathbb{Z})/\SL(2,\mathbb{Z})$ and boundary amenability}

Bekka and Kalantar discovered in \cite{bk} a deep relationship between properties of a subgroup $H$ of a discrete group $G$ and the ideal structure of the image of the full group $C^*$-algebra of $G$ under the corresponding quasi-regular representation $\lambda_{G/H}$. In particular they isolated two classes of subgroups (among others), namely the class $Sub_{sg}(G)$ of subgroups with the spectral gap property, for which $\lambda_{G/H} (C^*(G))$ contains a minimal (non-trivial, closed) ideal and the class $Sub_{w-par}(G)$ of weakly parabolic subgroups, for which $\lambda_{G/H}(C^*(G))$ contains a maximal ideal. In \cite{bk} Remark 7.6 the authors speculate about the existence of a group $G$ and a subgroup $H \in Sub_{sg}(G) \cap Sub_{w-par} (G)$ for which the two ideals coincide.\\
In the following we show that an example of such a situation is provided by the copy of $\SL(2,\mathbb{Z})$ inside  $\SL(3,\mathbb{Z})$ given by
\begin{equation*}
\left( \begin{array}{cc}	1	&	\begin{array}{cc}0	& 0 \end{array}\\
			\begin{array}{c}0\\ 0 \end{array}  &	\SL(2,\mathbb{Z})\end{array}\right).
			\end{equation*}
 
The proof relies on a dynamical presentation of $\SL(3,\mathbb{Z}) /\SL(2,\mathbb{Z})$ as an orbit of the diagonal action of $\SL(3,\mathbb{Z})$ on itself. The boundary of this orbit sits in a suitable boundary piece introduced in \cite{bip}, from which it follows that this action is boundary amenable. Hence the $C^*$-simplicity of $\SL(	3,\mathbb{Z})$ gives the desired example.\\
We retain the notation of \cite{bk}, \cite{bip} and refer to \cite{a}, \cite{bo} for the definition of amenable action and related topics.

  If $Z$ is a discrete countable space, we denote by $\Delta_\beta Z$ its Stone-{\v C}ech compactification and by $\partial_\beta Z = \Delta_\beta Z \backslash Z$ its Stone-{\v C}ech boundary. Recall that $\Delta_\beta Z$ is the set of ultrafilters on $Z$ with the topology generated by the basis of open sets of the form $U_E:=\{ \omega \in \Delta_\beta Z \; | \; E \in \omega\}$, for $E\subset Z$; the Stone-{\v C}ech boundary $\partial_\beta Z$ is the set of free ultrafilters. Every function $f$ from $Z$ to a compact Hausdorff space admits a unique continuous extension to the Stone-{\v C}ech compactification given by $f(\omega) = \lim_{z \rightarrow \omega} f(z)$. In particular group actions on $Z$ extend	 to $\Delta_\beta Z$.

The following is an application of the results contained in \cite{bk}. We recall that a group $G$ is $C^*$-simple if its reduced group $C^*$-algebra $C^*_r (G)$ is simple.
\begin{proposition}
\label{prop0}
Let $G$ be a discrete countable group and $H \in  {Sub}_{sg}(G)$. If $G$ is $C^*$-simple and the action of $G$ on the Stone-{\v C}ech boundary of $G/H$ is amenable, then $ \lambda_{G/H} (C^*(G))$ contains the compact operators as unique non-trivial closed ideal.
\end{proposition}
  \proof In virtue of \cite{bk} Theorem 7.4, $C^*(\lambda_{G/H})$ contains the compact operators as a minimal ideal $I_{min}$. Since the action of $G$ on $G/H$ is boundary amenable, the $*$-homomorphism $C^*(\lambda_{G/H} (G)) \rightarrow   C^*(\lambda_{G/H} (G))/I_{min}$ factors trough $C^*_r (G)$ (this observation is contained in the proof of \cite{Oz} Proposition 4.1 (1)$\Rightarrow$(3) and the proof of \cite{a} Proposition 5.13 (1)$\Rightarrow$(2)). Hence $C^*(\lambda_{G/H} (G))/I_{min} \simeq C^*_r (G)$ since $G$ is $C^*$-simple. In particular $C^*(\lambda_{G/H} (G))/I_{min}$ does not contain non-trivial ideals and so $I_{min}$ is the only ideal in $C^*(\lambda_{G/H} (G))$. $\Box$\\
  
  In order to apply the above result for the case at hand, we proceed in showing amenability of the action of $\SL(3,\mathbb{Z})$ on the Stone-{\v C}ech boundary of $\SL(3,\mathbb{Z})/\SL(2,\mathbb{Z})$. \\
  
  Let $\Gamma$ be a discrete countable group acting on a compact Hausdorff space $X$. We denote by $\mathcal{P}(\Gamma)$ the set of probability measures on $\Gamma$ endowed with the action of $\Gamma$ given by $g \mu = \mu \circ g^{-1}$. We recall from \cite{a} Definition 5.1 that an approximate invariant continuous mean for the action is a net $\mu_\lambda: X \rightarrow \mathcal{P}(\Gamma)$ satisfying $\lim_{\lambda} \sup_{x \in X}\| g\mu_\lambda (x)-\mu_\lambda (gx)\|_1=0$ for every $g \in \Gamma$, where $\|\cdot \|_1$ denotes the $l^1$-norm. The action is called amenable if it admits an approximate invariant continuous mean. If a set $Z$ is endowed with an action of $\Gamma \times \Gamma$, we will refer to the diagonal action of $\Gamma$ as to the action of $\Gamma$ on $Z$  given by $(\gamma, z) \mapsto (\gamma , \gamma) \cdot z$.
  \begin{lemma}
  \label{lem1}
  Let $\Gamma$ be a discrete countable group, $X$ a compact $\Gamma \times \Gamma$-space and $Y$ a compact $\Gamma$-space. Suppose the $\Gamma$-action on $Y$ is amenable and there is a continuous map $\phi : X \rightarrow Y$ such that $\phi ((g,h) x)=g\phi (x)$ for every $(g,h) \in \Gamma \times \Gamma$ and every $x \in X$. Then the diagonal action of $\Gamma$ on $X$ is amenable.
  \end{lemma}
  \proof Let $\mu_\lambda : Y \rightarrow \mathcal{P}(\Gamma)$ be an approximate invariant continuous mean. Then $\diag (\mu_\lambda \circ \phi)$ is an approximate invariant continuous mean for the diagonal action of $\Gamma$ on $X$, where $\diag : \mathcal{P} (\Gamma) \rightarrow \mathcal{P}(\diag (\Gamma))$ is given by $\diag (\mu)(g,g)=\mu (g)$. $\Box$\\

 \begin{theorem}
 \label{thm1}
The action of $\SL(3,\mathbb{Z})$ on the Stone-{\v C}ech boundary of $\SL(3,\mathbb{Z})/\SL(2,\mathbb{Z})$ is amenable.
\end{theorem}
\proof Consider the diagonal action of $\SL(3,\mathbb{Z})$ on itself, which is given by $(g, x) \mapsto  gxg^{-1}$ for $g,x \in \SL(3,\mathbb{Z})$. The stabilizer of the point
\begin{equation*}
z_0=\left(\begin{array}{ccc}
 1	&	0	&	0\\
0	&	-1	&	0\\
0	&	0	&	-1
\end{array}\right)
\end{equation*}
is a diagonal copy of $\SL(2,\mathbb{Z})$. Hence the orbit $Z$ of $z_0$ is a dynamical system isomorphic to $\SL(3,\mathbb{Z})/\SL(2,\mathbb{Z})$ and we can identify the Stone-{\v C}ech boundary of $\SL(3,\mathbb{Z})/\SL(2,\mathbb{Z})$ with the closed subset of $\partial_\beta \SL(3,\mathbb{Z})$ given by $\partial_\beta Z = \{ \omega \in \partial_\beta \SL(3,\mathbb{Z}) \; | \; F \cap Z \neq \emptyset \; \forall F \in \omega\}$. 
Let $M_3 (\mathbb{R})_{\leq 1}$ be the compact set of $3 \times 3$ real matrices of norm at most $1$ and $\phi : \SL(3,\mathbb{R}) \rightarrow M_3 (\mathbb{R})_{\leq 1}$ be given by $\phi (g)= g /\|g\|$. Every $z\in Z$ satisfies $z=z^{-1}$ (since $z_0^{-1} =z_0$), it follows that for every diverging sequence $\{z_n\} \subset Z$ for which $\phi (z_n)$ converges we have $(\lim_n \phi(z_n))^2=\lim_n (\phi (z_n))^2=0$. Let $\omega \in \partial_\beta (Z) \subset \partial_\beta \SL(3,\mathbb{Z})$. For every $g \in \SL(3,\mathbb{Z})$ let $s_{i }(g)$, $i=1,2,3$ be the eigenvalues of $\sqrt{g^t g }$, ordered in such a way that $s_{1 }(g) \geq s_{2 }(g) \geq s_{3 }(g)$ (cf. \cite{bip} 6.1). Consider the functions $  f_{1,2}: g \mapsto  s_1(g)/s_2(g) $ and $ f_{2,3} : g \mapsto  s_2(g) / s_3(g) $. We keep the same notation for their extensions to functions from $\Delta_\beta (\SL(3,\mathbb{Z}))$ to $[0,\infty]$. Suppose that $f_{1,2} (\omega) = c < \infty$. Then there is a sequence $(g_n)_n$ in $\SL(3,\mathbb{Z})$ with $\lim_n f_{1,2} (g_n) =c$ and $\lim_n \phi (g_n)^2=0$. Following \cite{bip} 6.1 for every $n \in \mathbb{N}$ let $s_{i,n} := s_i (g_n)$ for $i \in \{1,2,3\}$ and let $s_n$ be the diagonal matrix with entries $(s_n)_{i,i} = s_{i,n}$; we can write 
 \begin{equation*}
g_n=a_n s_n b_n = a_n\left( \begin{array}{ccc}	s_{1,n}	&	0	&	0	\\
							0		&	s_{2,n}	&	0	\\
								0	&	0		&	s_{3,n}\end{array}\right)b_n,
\end{equation*}
where $a_n, b_n \in \SO(3,\mathbb{R})$. Up to taking a subsequence, we can suppose that $a_n \ra a$, $b_n \ra b$ in $\SO(3,\mathbb{R})$, thus also $\lim_n \phi (s_n)$ exists and $$\lim_n \phi (a_n s_n b_n)= \lim_n a_n \phi (s_n) b_n= a\lim_n \phi(s_n) b;$$ the condition $\lim_n f_{1,2} (g_n) =c$ entails $$\rank (\phi (\omega))=\rank(\lim_n s_n/\|s_n\|))=2,$$ but then we should have $(\lim_n \phi(g_n))^2 \neq 0$, which is impossible. In the same way, if $ f_{2,3} (\omega) =c < \infty$, then there is a sequence $g_n$ in $\SL(3,\mathbb{Z})$ such that $\phi (\omega) = \lim_n \phi (g_n)$ and $f_{2,3} (\omega)=\lim_n f_{2,3} (g_n)$. But then $\lim_n \phi (g_n^{-1})$ has rank $2$, which again is impossible. Hence $\partial_\beta Z \subset \partial_0  \SL(3,\mathbb{Z}) $ in virtue of \cite{bip} Lemma 6.1, where $\partial_0 \SL(3,\mathbb{Z})$ is the boundary piece associate to the subgroup $P_0 \subset \SL(3,\mathbb{R})$ of upper triangular matrices. Since the action of $\SL(3,\mathbb{Z})$ on $\SL(3,\mathbb{R})/P_0$ is amenable, it follows from \cite{bip} Lemma 3.10 that the left action of $\SL(3,\mathbb{Z})$ on $C(\partial_0 \SL(3,\mathbb{Z}))^{\SL(3,\mathbb{Z})_r}$ is amenable; it follows from Lemma \ref{lem1} that the diagonal action of $\SL(3,\mathbb{Z})$ on $\partial_0 \SL(3,\mathbb{Z})$ is amenable (the action of $\SL(3,\mathbb{Z}) \times \SL(3,\mathbb{Z})$ on both $\partial_0 \SL(3,\mathbb{Z})$ and $\partial_\beta \SL(3,\mathbb{Z})$ is given by $(g_1, g_2) \cdot \omega = g_1 \omega g_2^{-1}$). Hence the diagonal action of $\SL(3,\mathbb{Z})$ on $\partial_\beta Z$ is amenable. $\Box$

\begin{corollary}
The $C^*$-algebra $ \lambda_{\SL(3,\mathbb{Z})/\SL(2,\mathbb{Z})} (C^*(\SL(3,\mathbb{Z}))$ contains the compact operators as unique non-trivial closed ideal.
\end{corollary}
\proof In \cite{bk} Example 7.7 it is shown that $\SL(2,\mathbb{Z})  \in Sub_{sg}(\SL(3,\mathbb{Z})) \cap Sub_{w-par}(\SL(3,\mathbb{Z}))$. The $C^*$-simplicity of $\SL(3,\mathbb{Z})$ is proved in \cite{bch}. Hence the result follows from an application of Theorem \ref{thm1} and Proposition \ref{prop0}. $\Box$
\bibliographystyle{mscplain}
 \bibliography{biblio}
 
\section{Acknowledgments}

The authors thank the anonymous referee for the constructive comments on a previous version of the article. The authors acknowledge the support of INdAM-GNAMPA and the MIUR Excellence Project awarded to the Department of Mathematics, University of Rome $2$, CUP E83C180001000 and of the grant Beyond Borders CUP: E89C20000690005. The first named author was supported by the MIUR grant CUP: E83C18000100006 and by the grant Beyond Borders CUP: E84I19002200005. The second named author acknowledges the partial  support of the grant CNCS Romania, PN-III-P1-1.1-TE-2019-0262. The present project is part of: - OAAMP - CUP E81I18000070005.

\baselineskip0pt

\end{document}